\documentclass[twoside,leqno]{article}

\usepackage[letterpaper]{geometry}

\usepackage{siamproceedings}

\usepackage[T1]{fontenc}
\usepackage{amsfonts}
\usepackage{graphicx}
\usepackage{epstopdf}
\usepackage{enumitem}
\usepackage{algorithmic}
\ifpdf
  \DeclareGraphicsExtensions{.eps,.pdf,.png,.jpg}
\else
  \DeclareGraphicsExtensions{.eps}
\fi

\newsiamremark{remark}{Remark}
\newsiamremark{hypothesis}{Hypothesis}
\crefname{hypothesis}{Hypothesis}{Hypotheses}
\newsiamthm{claim}{Claim}
\newsiamremark{example}{Example}
\newsiamremark{question}{Question}

\usepackage{amsopn}

\begin{document}

\newcommand\relatedversion{}
\renewcommand\relatedversion{\thanks{The full version of the paper can be accessed at \protect\url{https://arxiv.org/abs/0000.00000}}} 

\title{\Large AI and humanity's long conversation}
    \author{Geordie Williamson\thanks{Sydney Mathematical Research Institute (\email{g.williamson@sydney.edu.au}, \url{https://www.maths.usyd.edu.au/u/geordie/}).}
}

\date{}

\maketitle


\fancyfoot[R]{\scriptsize{Copyright \textcopyright\ 20XX by SIAM\\
Unauthorized reproduction of this article is prohibited}}





\begin{abstract} These are notes to accompany my public lecture at the International
  Congress of Mathematicians, in Philadelphia in July, 2026. The goal
  is to discuss some developments at the interface of AI and mathematics in the context of
  the long conversation of mathematics.
\end{abstract}

\section{Introduction.}
When I was asked to give a talk at the interface of AI and mathematics, the first quote
that came to mind is from Mazur \cite{mazur}:
\begin{quote}
“...there is an underlying timelessness in the basic conversation that is mathematics.” 
\end{quote}
I am not sure who first came up with the idea of mathematics as a
``long conversation'', but it articulates one of the things that I love most about
the subject. I love the picture of mathematics as a conversation
stretching over millennia. I love the flurry of activity when a new approach
emerges. I love the long silences when we don't know what to say, or
are digesting what has been said by others. I love the slow digestion of the
salient points into our mathematical culture. Like many others, I get
great motivation from the possibility of contributing something
interesting to this long conversation.

This lecture concerns the interaction of this patient and slow
conversation with the current breathless conversation around AI.
Progress in AI in mathematics over the last few months has been both surprising and
confronting. It is
incredible for me to think that when I began preparing this lecture in early May
(only four months ago), I felt the need to prepare several slides to
try to convince the working mathematician that LLMs could sometimes be
useful tools in their research!
It is hard not to see current advances in AI as a turning point in this long conversation, and it is hard to predict how
mathematics will
adapt. 

The goal of this lecture is not to make predictions, but to point out
things that are already evident (and were already
evident a year or two ago). This lecture is structured around two aspects:

\subsection{The wonder.} 
The first is a sense of \emph{wonder} at this technology. I have felt
for several years now that mathematicians should be more aware of
neural networks and developments in AI, both as sources of
interesting problems and as tools for research. Of course, recent
developments have forced us to pay attention. However, it is easy to forget that
the foundations of modern AI are mathematical, and involve beautiful
ideas that are comparatively easy for the mathematician to
comprehend. There are many mysteries surrounding neural networks
leading to interesting research directions. In
the first section of these notes I try to explain the basic ideas
underlying LLMs with the working mathematician in mind.

\subsection{The concern.} 
The second is a deep \emph{concern} about this technology. It has
already changed the working lives of many software engineers,
and looks to be changing our lives as mathematicians too. The last
months have seen many discussions of what the mathematician does, and
what they are likely to do in the future. Here I would like to
emphasise three aspects:\footnote{I promise these three
points were not as obvious in July as they appear now.}
\begin{enumerate}
\item We should expect significant change, and that it will not be easy for the
  mathematical community to adapt;
\item The best
  guess is that models will be more capable in the months and
  years to come;
\item These developments force us to honestly confront what we do
  and value as mathematicians.
\end{enumerate}
Over the last two months, the \emph{concern} has (justifiably!) come to dominate the conversation, and many
accounts have been written along these lines. The reader interested in
discussions about the mathematical process, and how it should respond
to AI is referred (for example) to the excellent essays \cite{Bessis,
  Gamburd,TaoAI}, as well as to the many contributed posts on
\emph{Proofs and Prompts} \cite{ProofsAndPrompts} and Terry Tao's blog 
\cite{TaoBlog}.

In the following, I focus more on the \emph{wonder} and discuss the \emph{concern} more briefly, in the same terms it
was presented in July.

\section{The wonder: next token prediction and transformers.}

In 2026, most mathematicians have had the experience of putting a query
into an LLM and getting a useful answer. It is surprising
that this is possible, and many (myself included) predicted as recently
as two years ago that it would not be. Here I try to outline some of
the mathematical ideas underlying LLMs, trying to
emphasise both the simplicity of the underlying ideas and the
persistent wonder that they work at all.

\subsection{Shannon's series of approximations to English.}

Let us begin with Shannon's remarkable paper \emph{A mathematical
  theory of communication}, 
published in 1948 \cite{shannon}.\footnote{Incidentally, this was the same year
that Turing wrote \emph{Intelligent Machinery} \cite{Turing1948} which
is one of the first articles that takes the prospect of thinking
machines seriously. In this article he identifies mathematics as an
interesting testing ground for such ideas, which turned out to be a
very prescient prediction.}
Shannon's concern is how efficiently information may
be transmitted. The more structured
information is, the more efficiently it can be compressed and hence
transmitted. In Shannon's words, increased efficiency is possible ``\dots due to the statistical
structure of the original message''.

\begin{figure}[tbp]
  \centering
  \includegraphics[width=0.8\columnwidth]{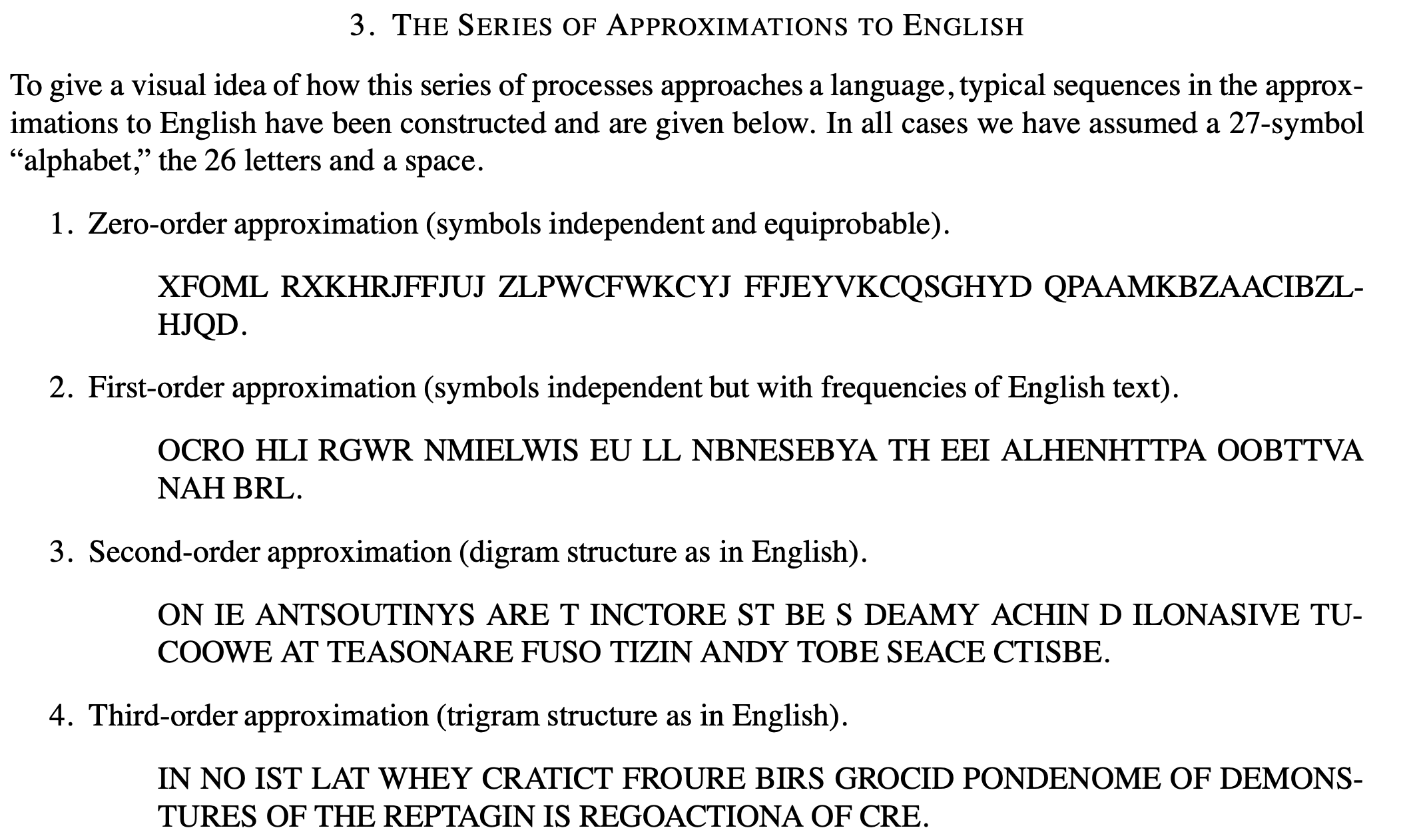}
  \caption{Shannon's ``series of approximations to English'', from \cite[\S~3]{shannon}.}
  \label{fig:shannonapprox}
\end{figure}

He first considers various simple models
(e.g. uniform distributions and Markov chains) for which the regularity is manifest in the
definition. He then asks himself to what extent these simple models
describe the structure present in English. Based on a simple
model consisting of 26 letters and a space, he considers
how well simple models like that of a Markov chain approximate English
(see \cref{fig:shannonapprox}). In what Shannon terms the ``zero order'' approximation, each
symbol is sampled uniformly at random. In the ``first order''
approximation, each symbol is sampled according to its frequency in
English. In the ``second order'' approximation, each symbol is sampled
according to the frequency with which it follows the previously
sampled symbol (a \emph{bigram} distribution).

He then moves on to a series of approximations with letters replaced by
words. Here is his sample from a bigram distribution on words:
  \begin{quote}
    THE HEAD AND IN FRONTAL ATTACK ON AN ENGLISH WRITER THAT THE
    CHARACTER OF THIS POINT IS THEREFORE ANOTHER METHOD FOR THE
    LETTERS THAT THE TIME OF WHO EVER TOLD THE PROBLEM FOR AN
    UNEXPECTED.
\end{quote}
In order to produce the bigram sample, he uses the following simple
and clever sampling
technique. He first chooses a book at random from the bookshelf, and
opens the book at a random page. He then reads until he finds the
word in question. He then writes down the next word. This gives a
good approximation of the bigram distribution in English.

It is interesting to reflect that, even if Shannon were to use the
internet as his bookshelf, even short and reasonable sequences of
words simply do not occur. So a naive approach to $n$-gram
estimation for small $n$ is impossible, even with the internet.
For example, when this article was written,
Google attests that the three words ``tickled an echidna'' do
not occur on the (Google-)known internet! It is easy to come up with other
examples. Thus in order to approximate the distribution of English we
need more sophisticated techniques.

\begin{example}
  Despite its difficulties, $n$-gram approximation has been a powerful
  tool in older language models used for machine translation. For example,
  \cite{brants2006web1t} contains the $5$-gram
  frequencies for over one billion sequences found on the internet in 2006.
\end{example}

\subsection{Next token prediction.}

In modern language, what Shannon is attempting to do here is
\emph{next token prediction}. We first fix a finite
set of tokens, which are typically chunks of words, but also
incorporate punctuation, start and end tokens etc.:
\[
\text{`` tick '', `` ech '', ``id '', ``'ve '', ``! '' , ``<END>''},
\dots
\]
Our tokens are chosen so that we can break up English sentences efficiently into
sequences of tokens:
\[
\text{``I've never tickled an echidna!''} \to
[\text{``I'', ``'ve'', `` never '', `` tick'', ``led'', `` an '', `` ech'', ``id'',
  ``na'', ``!''} ]
\]
The task of next token prediction is to predict the probability of the
next token, given all that has been read so far.

\begin{example}
  As a first example, consider the delightful opening to
\emph{Philosophy and Fun of Algebra}, written by Mary Everest Boole in
1909 \cite{Boole}. Your task is to predict the word (``token'') hiding behind the box:
\begin{equation*}
  \includegraphics[width=12cm]{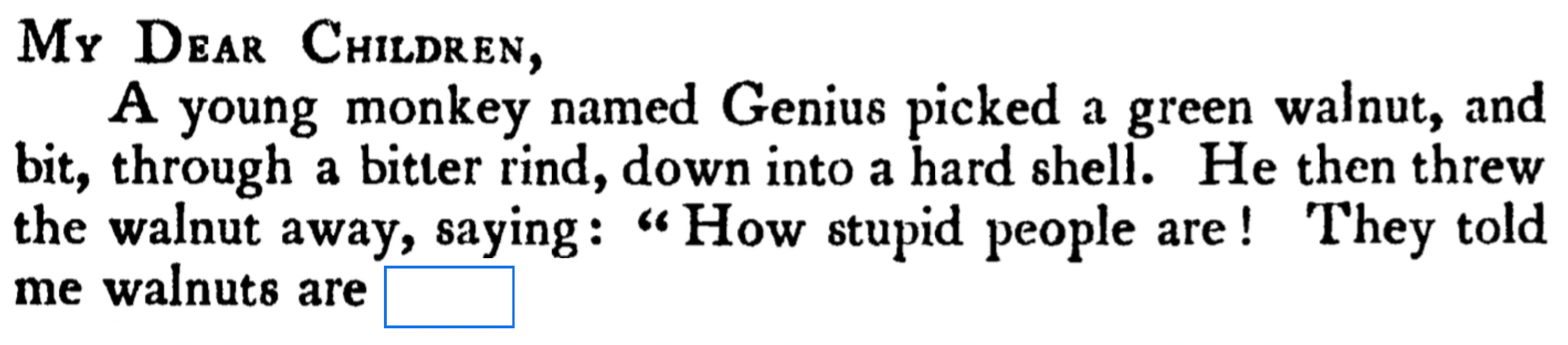}
\end{equation*}
In this example, you have very little context. You know little about
this monkey called Genius, and though you can make some reasonable
guesses about what follows, you are probably far from sure.
This example illustrates what most people imagine when they
first think about next token prediction. It is a sort of statistical
game, and something implementing it is a ``stochastic
parrot''. \end{example}

\begin{example} Our second example is from a book that caused me
much suffering as a graduate student, Hartshorne's
\emph{Algebraic geometry} \cite{Hartshorne}. Again your task is to predict the
next token:
\begin{equation*}
  \includegraphics[width=12cm]{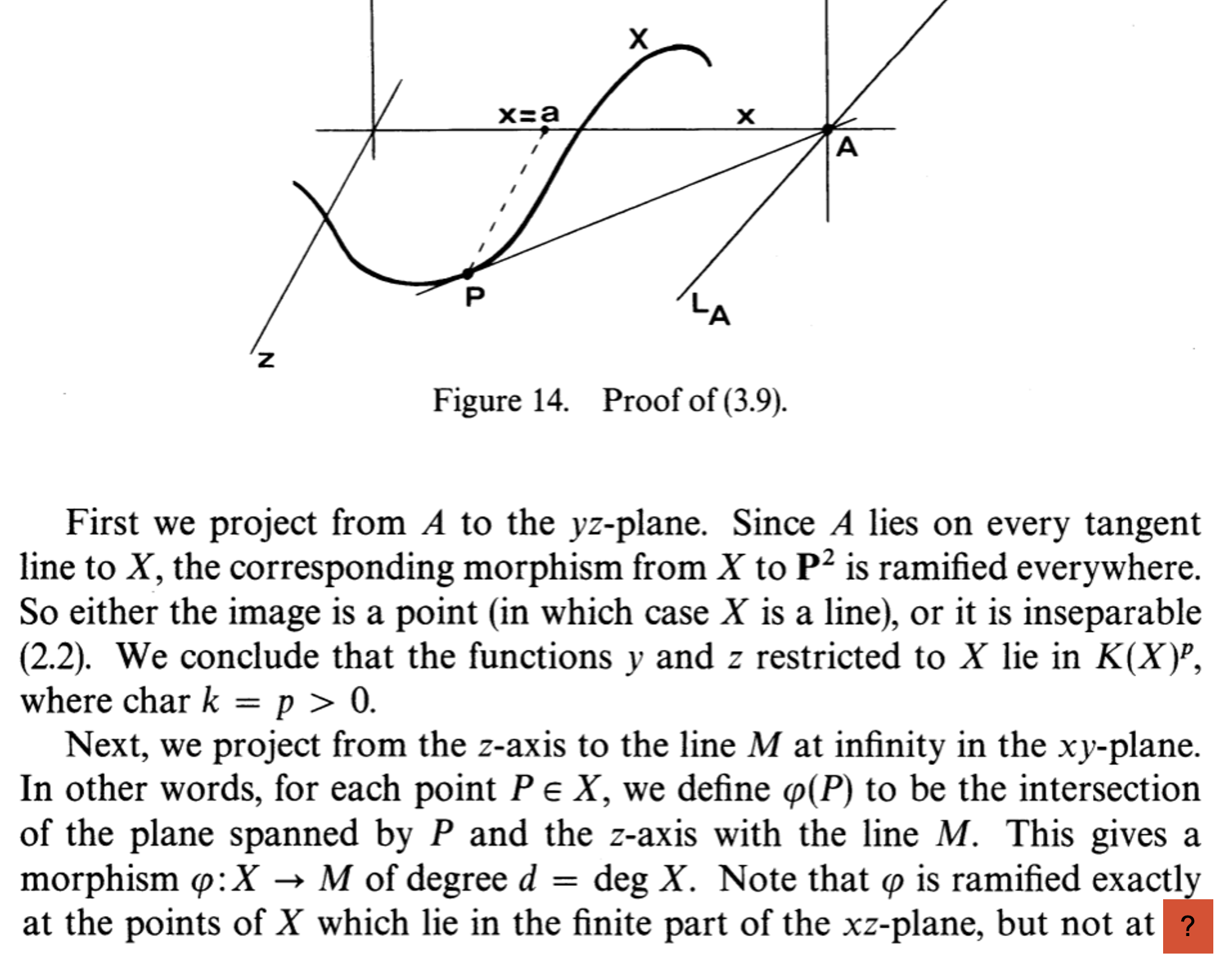}
\end{equation*}
If you are a mathematician, you can guess that this looks like a
proof and the missing token is almost certainly a symbol indicating a
point. However, even with the context of the paragraph containing the
missing token we cannot be sure what point is meant. However, if you
know some geometry and can read the previous page, you can be nearly
100\% sure of what the missing token is.

This example is meant to challenge the ``stochastic parrot'' point of
view. One can encode arbitrary amounts
of understanding in next token prediction. If a system is very good at
next token prediction, it necessarily has some sort of understanding.
\end{example} 

\subsection{Transformers.} We have already attempted to point out the
difficulty of constructing further approximations to English via
trigrams, quadgrams etc. Shannon finishes this section of his 1948 paper
with the prognosis: ``It would be interesting if further approximations
could be constructed, but the labor involved becomes enormous at the
next stage.'' Ever since Shannon's work, it has been a major challenge
to construct better approximations to language (see e.g. \cite{Goodman2001Bit}).

\begin{figure}[tbp]
  \centering
  \includegraphics[width=0.49\columnwidth]{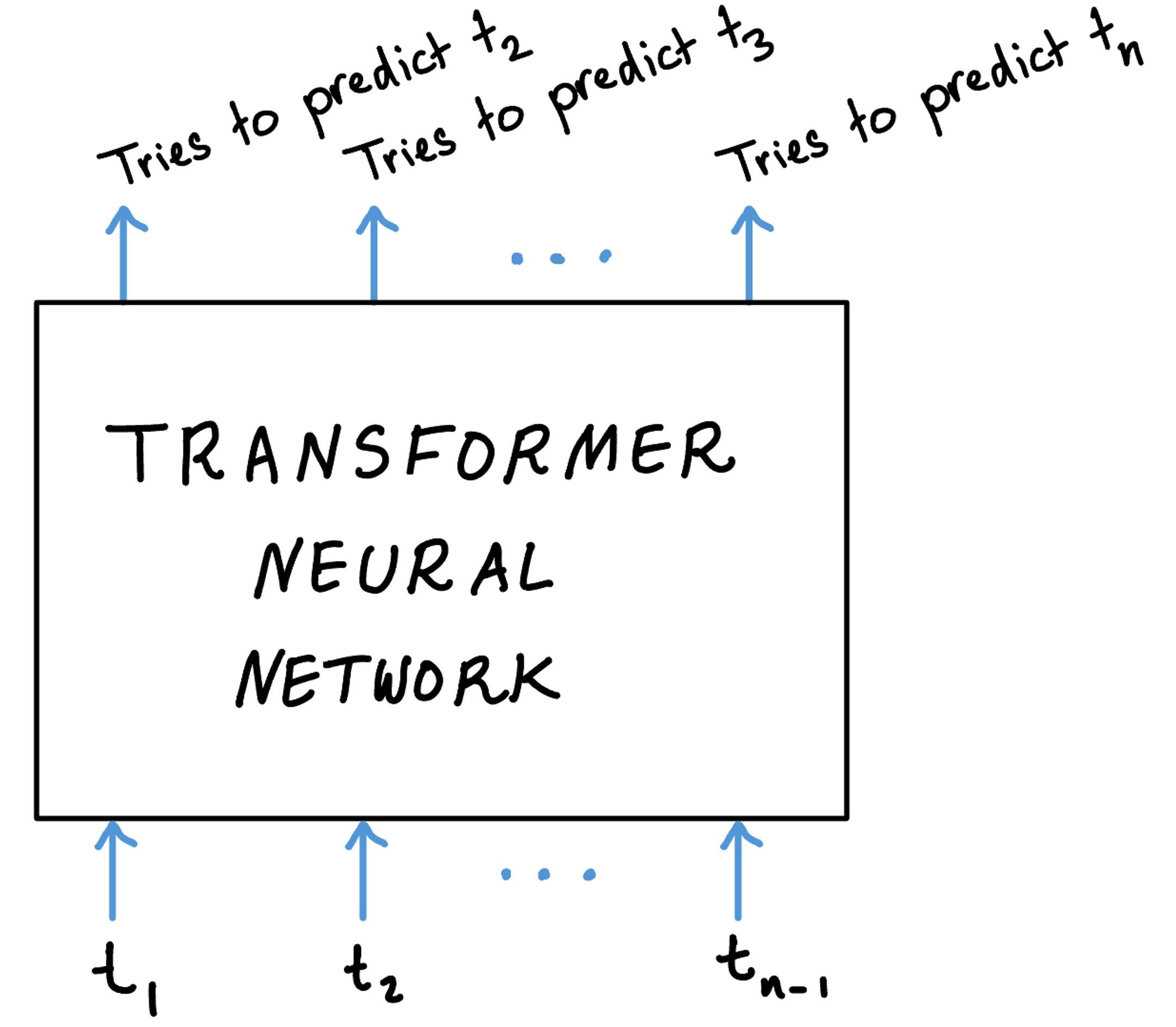}
  \caption{The simplest mental model of a transformer.}
  \label{fig:transformer}
\end{figure}

This is exactly the problem on which transformers have excelled. A
\emph{transformer} is a type of neural network. A basic mental model
for a transformer is displayed in \cref{fig:transformer}. The input is a (variable length) sequence of tokens
\[
(t_1, t_2, \dots, t_{n})
\]
and the output is a sequence of probability distributions
\[
(p_1,p_2,\dots, p_n)
\]
over the set of tokens. These distributions are meant to predict what
probability the model assigns to the possible next tokens. For example, if our vocabulary has three
tokens then $p_1$ might be $(0.01,0.1,0.89)$, indicating that the
model thinks that the first token is very unlikely to follow $t_1$, the 
second token unlikely, and the third token very likely. Note that the
transformer is making many predictions (i.e. $p_1$, $p_2$, \dots) of
possible subsequent tokens in parallel.\footnote{Note however that
  every prediction is conditioned on the preceding tokens. This means
  that during training many predictions are made in parallel, but when
  sampling one needs to perform $n$ sequential forward passes through
  the transformer. This explains why input tokens are typically much
  cheaper than output tokens.}

\begin{figure}[tbp]
  \centering
  \includegraphics[width=0.49\columnwidth]{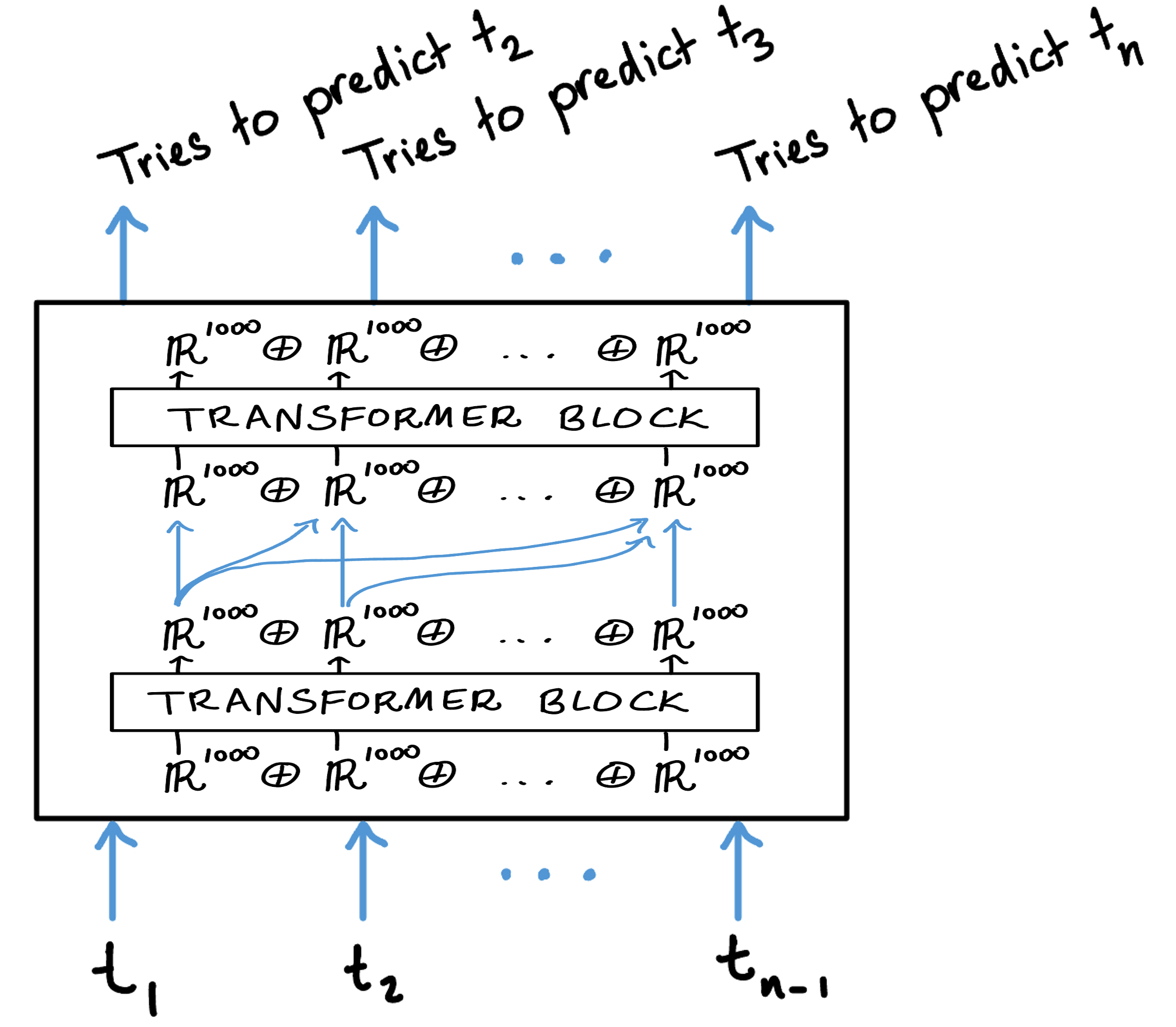}
  \caption{A more accurate mental model of a transformer.}
  \label{fig:transformer2}
\end{figure}

So what is inside the transformer? A big pile of linear algebra! I
don't want to go into too much detail here, but a good mental model to
have in mind is that at the first stage the sequence of tokens is
transformed into a sequence
\[
(v^0_1, ..., v^0_{n}) \in \mathbb{R}^{N} \oplus \dots \oplus \mathbb{R}^{N} 
\]
of $n$ vectors inside some fixed high-dimensional vector
space $\mathbb{R}^{N}$ (the \emph{embedding space}). These sequences are then processed several times by
\emph{transformer blocks}. Passing through each transformer block
transforms a sequence
\[
  (v^{i-1}_1, ..., v^{i-1}_{n}) \in \mathbb{R}^{N} \oplus \dots \oplus \mathbb{R}^{N} 
\]
into a new sequence of vectors
\[
(v^i_1, ..., v^i_{n}) \in \mathbb{R}^{N} \oplus \dots \oplus \mathbb{R}^{N}.
\]
Finally, after passage through some number of transformer blocks, this sequence of vectors is
decoded into probabilities. Each of the updates in the transformer
blocks is a relatively simple neural network, but involves many 
parameters. These parameters describe the affine linear maps which
take place within the transformer blocks, as well as learned
normalizations (see \cref{fig:transformer2}).\footnote{As in most neural
network, fixed non-linearities are also present.}

\begin{remark}
We will not give the gory details of transformer internals
here. For a long time, the canonical source was the original paper
\cite{attention} which is quite impenetrable for the mathematician without
substantial background in machine learning. However, recently several
other sources have appeared which have finally solved this exposition
problem. In my opinion, the best source is Sanderson's series of
videos
\cite{Sanderson2024Transformers,Sanderson2024Attention,Sanderson2024Facts}. The
reader wanting a detailed treatment in paper form is referred to
\cite{PhuongHutter}. For an interesting mathematical treatment which tries to
abstract away some of the engineering details, see
\cite{Mathematicaltrans}.
\end{remark}

It is best to think about a transformer (represented in
\cref{fig:transformer} and
\cref{fig:transformer2}) as representing a family of sequence
predictors that are controlled by many parameters.\footnote{Note that
  there is an ambiguity in the term transformer here. The same
  ambiguity is already present in the term \emph{neural network}. Typically, one
  uses \emph{neural network} to refer both to the family of
  programs or networks determined by parameters ``the neural network
  is initialized and then trained via gradient descent'' and a
  particular instantiation with parameters fixed ``the probabilities
  are calculated by passing an image through the neural network''.}
Typically, the transformer is initialized randomly. At initialisation
its output is close to a uniform distribution on sequences of
tokens.\footnote{This is somewhat problematic, as multiplicative
  discrepancies can accumulate on long sequences, but is a good first mental model.}
It
is then \emph{trained} on many samples from some distribution
(e.g. language from the internet) that one wants to approximate.

\subsection{Training transformers.} We now give a little more detail
on what the training procedure entails. Suppose that we have fixed a
transformer architecture which depends on parameters $\theta$. Let us
denote by
\[
p_{\theta}(t_i | t_1t_2 \dots t_{i-1})
\]
the probability that our model predicts $t_i$ to follow a sequence of
tokens  $(t_1, \dots, t_{i-1})$. Note that this is exactly one of the
outputs of our transformer. By the product rule for
probabilities, the probability that a sequence $(t_1, \dots, t_n)$
is sampled is given by:
\begin{equation} \label{eq:prob}\prod_{i=1}^{n} p_{\theta}(t_i | t_1t_2 \dots t_{i-1}).
  \end{equation}

Now suppose we are trying to approximate some distribution of tokens (e.g. English on
the internet). It is natural to try to adjust parameters to make samples from our
distribution slightly more likely. Because all parts of our
transformer are differentiable, we can ``bump'' our
transformer to slightly increase the probability of producing our
sample, via gradient ascent. Technically, it is more convenient to
consider the negative log of \eqref{eq:prob},
\begin{equation}
  \label{eq:cross-entropy}
- \sum_{i=1}^{n} \log p_{\theta}(t_i | t_1t_2 \dots t_{i-1})
\end{equation}
and to average it over a batch of samples from our distribution. This
yields the \emph{loss}, which we then try to reduce by many steps of
gradient descent.

Let us suppose for a moment that we succeed, and our transformer
perfectly matches the distribution from which our samples arise. In this case
the expected value of the loss (for sequences of length $n$) will be
\begin{equation}
  \label{eq:entropy}
- \sum_{(t_1,\dots,t_n)} p_{\theta}(t_1t_2 \dots t_n) \sum_{i=1}^{n} \log p_{\theta}(t_i | t_1t_2 \dots t_{i-1}).
\end{equation}
This is exactly the Shannon entropy of the distribution of sequences
of length $n$. More generally, the value of \eqref{eq:cross-entropy}
approximates the \emph{cross entropy} between the distribution
produced by the transformer, and the distribution one is attempting to
approximate. One can think of it as a measurement of the distance
between the distribution produced by our transformer, and the
distribution we are attempting to approximate.

\begin{example}
 It is fun and instructive to train a small transformer on a simple task. As
 it is training, one can sample from it and in some examples one gets
 a sense for its slow convergence to the target distribution. For
 example, following Karpathy \cite{KarpathyNanoGPT}, Papadopoulos \cite{PapadopoulosLLM} explains how a
 small transformer might be trained on Shakespeare.

 The transformer is
 then prompted with ``to be or not''. Here is a sample after 200 steps of gradient descent (the output
 looks essentially random):
 \begin{quote}
   To be or notis p ce mei odorethleedetire'ilethed ye m arkesothir fnon b tigb'i.
\end{quote}
Here is a sample after 1000 steps (some structure is emerging):
\begin{quote}
  To be or nothing are good men, \\
The profent of little, our actory.\\
CORIOLANUS:\\
Is it now of your many death?
\end{quote}
Finally, a sample at step 2400 (it is passably approximating Shakespeare):
\begin{quote}
  To be or not to be some of you shall know \\
That everlature by Romeo: what news,\\
Which you had knock'd my part to speak
\end{quote}
\end{example}

\subsection{Scaling laws.} If one experiments with transformers on
language, one notices consistent improvements with the size of the
transformer (number of parameters), amount of data, and the amount of
computation used for training (``compute''). It is then
natural to ask how well a large transformer trained for a long time on a large
chunk of the internet approximates English (or some other language, or
indeed all languages at once). Because
training a large transformer uses significant resources, it is natural
to ask whether the answer is predictable in advance.

\begin{figure}[tbp]
  \centering
  \includegraphics[width=1\columnwidth]{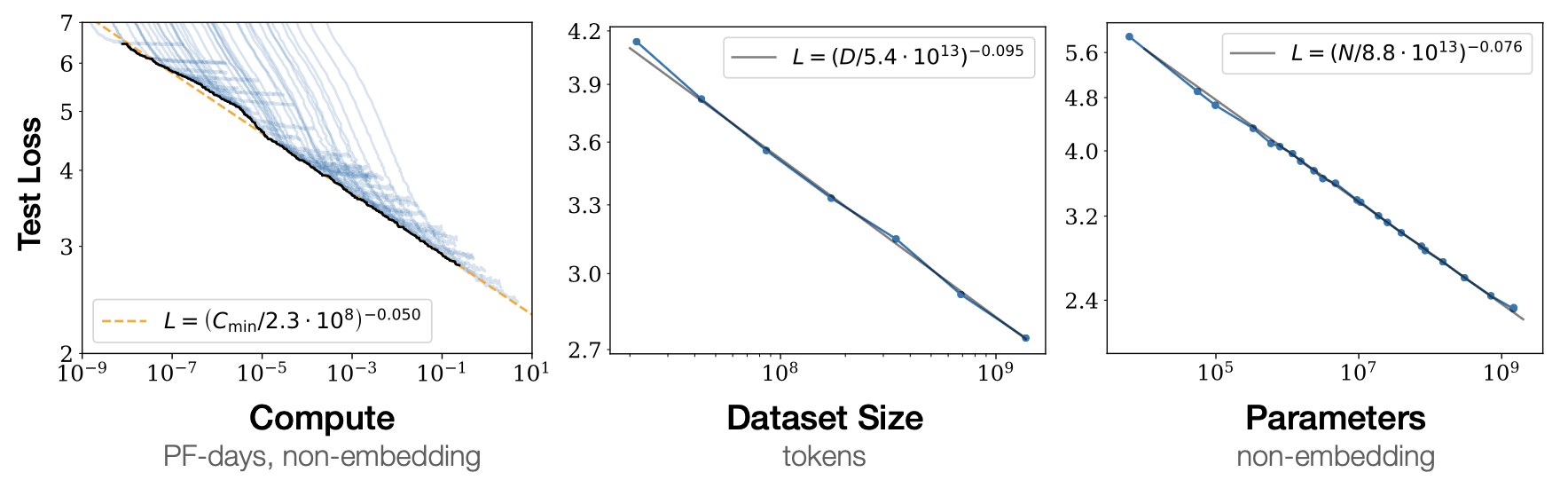}
  \caption{The famous scaling law plot, from the introduction to \cite{kaplan2020scaling}.}
  \label{fig:powerlaw}
\end{figure}

In order to measure this, it is natural to compute the cross-entropy
with the target distribution (e.g. samples from the internet). This
is approximated well by the \emph{loss} that was discussed in the
previous section.
This cross-entropy gives a concrete measurement of how
``far'' the model is away from exactly modelling the distribution, with
lower bound the entropy of the target distribution. The
remarkable empirical observation \cite{hestness2017deep,
  kaplan2020scaling, hoffmann2022training} is that this cross-entropy
obeys a power-law relationship in parameters, data and computational resources (see
\cref{fig:powerlaw}). These are the famous \emph{scaling laws}. In other words, by increasing parameters,
data and compute one obtains a predictable improvement in the ``distance''
to language.

Scaling laws explain to a large extent the
eye-watering investment that has taken place over the last five years
in AI and data-centres. Indeed, scaling laws give confidence in
a predictable increase in model quality.

\begin{remark}
  Nowadays, there is reliable empirical evidence for
  some form of scaling laws across many different aspects of machine
  learning. Typically, scaling laws take the form
  \[
L(N,D) = L_\infty + AN^{-\alpha} + BD^{-\beta}
\]
where $L(N,D)$ is the loss after training, $L_\infty$ is the entropy
of the target distribution, $N$ is number of parameters, $D$ is the
amount of training data, and $A$, $B$, $\alpha$ and $\beta$ are
constants. (Note that compute is typically proportional to $ND$.)
On the other hand there appears to be little theoretical understanding
of the origin of these scaling  laws, or the meaning of the exponents $\alpha$ and $\beta$ (see \cite[\S
2.3]{simon2026scientific}).
For a detailed discussion of scaling laws, see \cite{weng2026scaling}.
\end{remark}

\subsection{Post-training.} In LLM parlance, everything that we have
discussed so far belongs to \emph{pre-training}. A raw pre-trained LLM is a strange
beast. A useful, if imperfect, mental model is a slightly lossy
compression of the internet: it is likely to insult, lie or descend
into gibberish with the same probability of that taking place on the
part of the internet on which it was trained.

In order to convert a pre-trained LLM into something useful there is a
wide array of techniques called \emph{post-training}. These include
instruction tuning, preference optimisation, reinforcement learning, and
various forms of fine-tuning; many details of the post-training pipelines used
by leading laboratories appear to be closely guarded secrets. We won't
have much to say about that here, except to highlight that it is
post-training that converts a raw next-token predictor into something
that can write code which often compiles, and mathematics that is
sometimes correct.

\subsection{Many mathematical mysteries.} 
\begin{figure}[tbp]
  \centering
  \includegraphics[width=0.8\columnwidth]{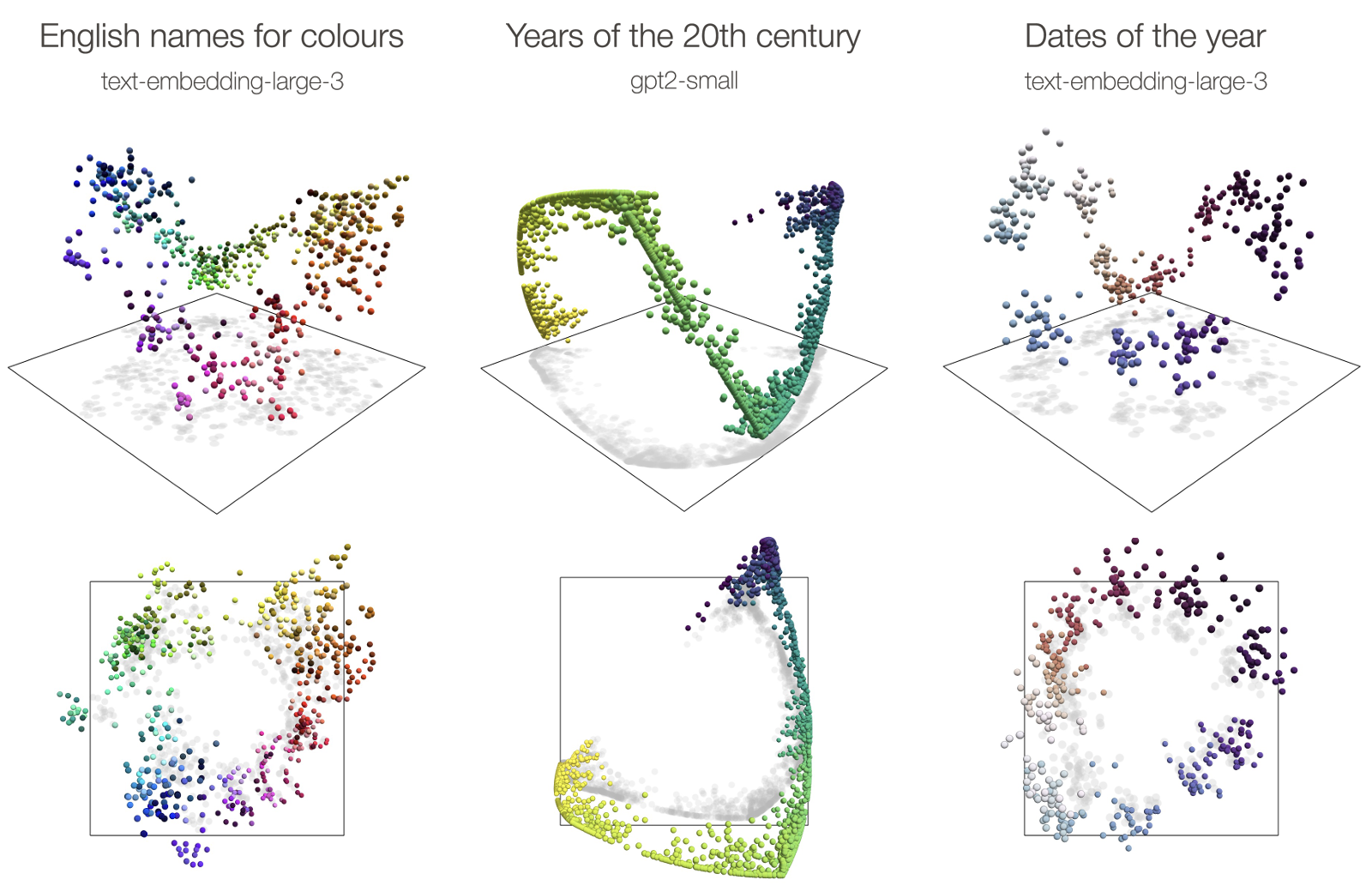}
  \caption{Emergent structure in embedding spaces of language models,
    from \cite{modell}.}
  \label{fig:embeddings}
\end{figure}

Although the basic ideas involved in modern LLMs are simple, there are
many mysteries surrounding generalization, internals and training dynamics. This is the subject of machine learning theory or
``math for AI'', which has blossomed into a flourishing subject.
Here are some examples of the interesting questions that arise:
\begin{enumerate}
\item Often neural networks are ``massively
  overparametrized''---the expressive power of neural
  networks far exceeds the amount of data on which the model is
  trained. Traditional statistical wisdom suggests that such
  models should overfit and generalize poorly (the ``bias/variance
  tradeoff'' -- think of fitting a high degree polynomial to a low number
  of points). The fact that neural networks often generalize well
  \emph{despite} overfitting has revealed a new
  paradigm in statistics, often discussed under the related notions of
  ``benign overfitting'' and ``double descent''
  \cite{belkin2019reconciling, nakkiran2020deep, bartlett2020benign,
    hastie2022surprises,belkin2021fit}.
  \item What is going on inside neural networks and transformers?
    There is now overwhelming empirical evidence that neural networks learn
    geometrically meaningful representations. For example, in
    \cite{nanda2023progress} the authors 
    establish that a transformer ``discovers the Fourier transform''
    in order to add integers modulo $113$. There is also evidence
    for similar structure emerging in much larger language models (see
    \cref{fig:embeddings}). The subject of \emph{mechanistic
      interpretability} attempts to decompose a trained neural network
    or transformer into small interpretable pieces, somewhat like decompiling
    machine code, or understanding brain function at the level of
    neurons. Although the subject is far from understanding what 
    is taking place in a large transformer, researchers have identified 
    interesting fundamental patterns and emergent phenomena \cite{olah2020zoom,
      elhage2021mathematical, olsson2022context}.
\item Why do scaling laws hold? Although there are now good theoretical
  models \cite{sharma2022scaling, bahri2024explaining,
    maloney2022solvable,hutter2021learning}, many mysteries
  remain. Despite the fact that scaling laws are reliable predictors
  of the loss of the trained model, the question of capabilities
  (e.g. what type of questions the model can correctly answer)
  remains subtle \cite{wei2022emergent,
    srivastava2023beyond,arora2023theory}.
\end{enumerate}

\begin{remark}
For an excellent survey of recent progress in theoretical and
empirical understanding of deep learning systems, the reader is
referred to \cite{simon2026scientific}.
\end{remark}

\begin{remark}
  The mysteries discussed above are not just of mathematical or
  scientific interest. Indeed, it is a priori clear
  that broadly capable AI systems carry a significant safety risk, and
  events this year have confirmed these dangers are real. Thus, better
  understanding of AI systems could lead to safer systems, and
  mathematicians can contribute to this research direction
  \cite{sharkey2025open, levine2026mathaisafety}.
\end{remark}

\subsection{LLMs and mathematics.} The last year has seen many
examples of LLMs contributing to mathematics, either (largely)
autonomously or through interaction with a mathematician. Here I would
like to discuss an early example of an interesting mathematical structure
discovered with the help of an LLM. I hope this simple
example will help the working mathematician gain some understanding of
``how'' LLMs are currently solving more significant problems. (Though the
methodology described here is different to most modern LLM
applications, it is much more concrete, and one can explicitly see
what input the LLM is providing.)

\begin{figure}[tbp]
  \centering
  \includegraphics[width=0.3\columnwidth]{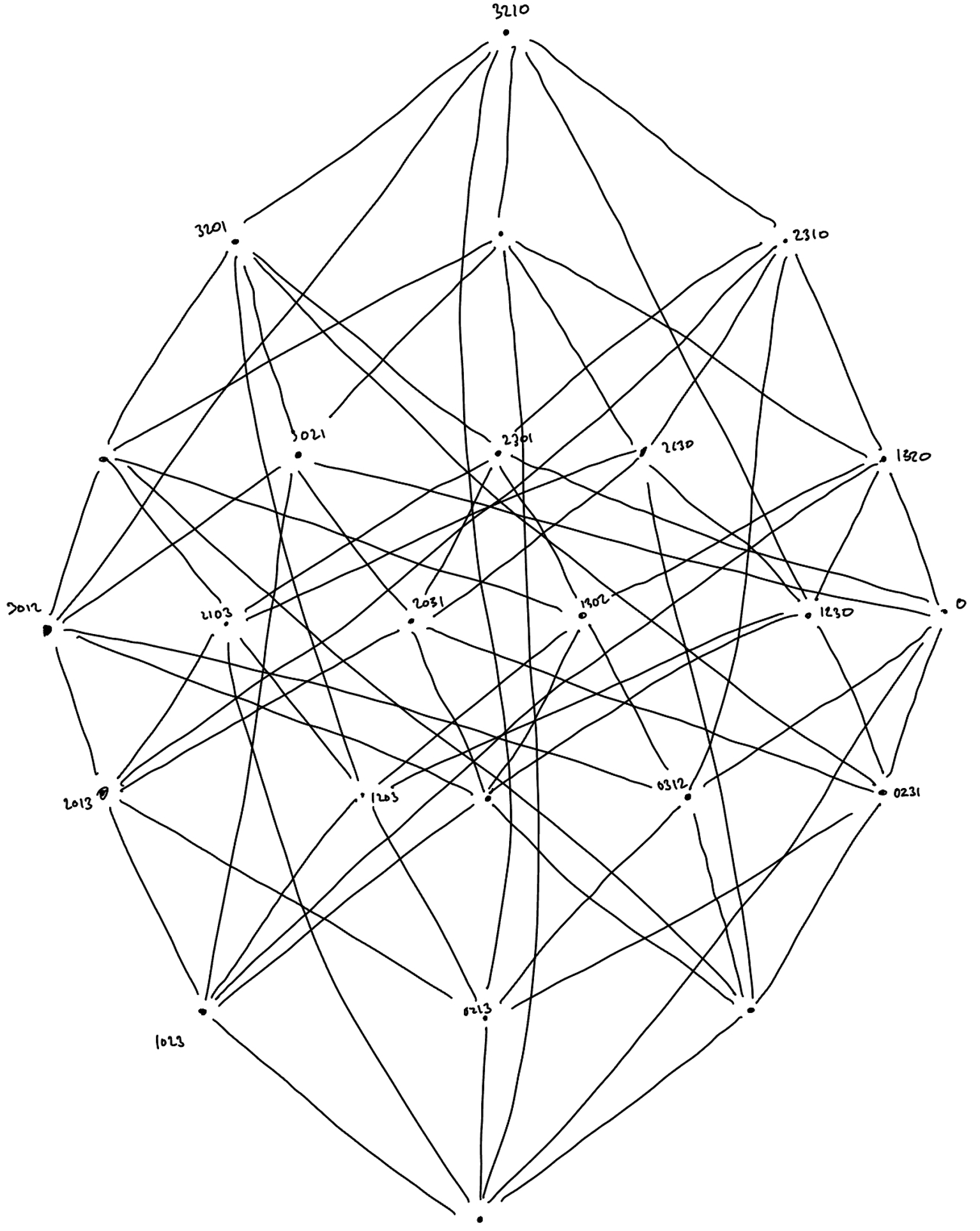}
  \includegraphics[width=0.3\columnwidth]{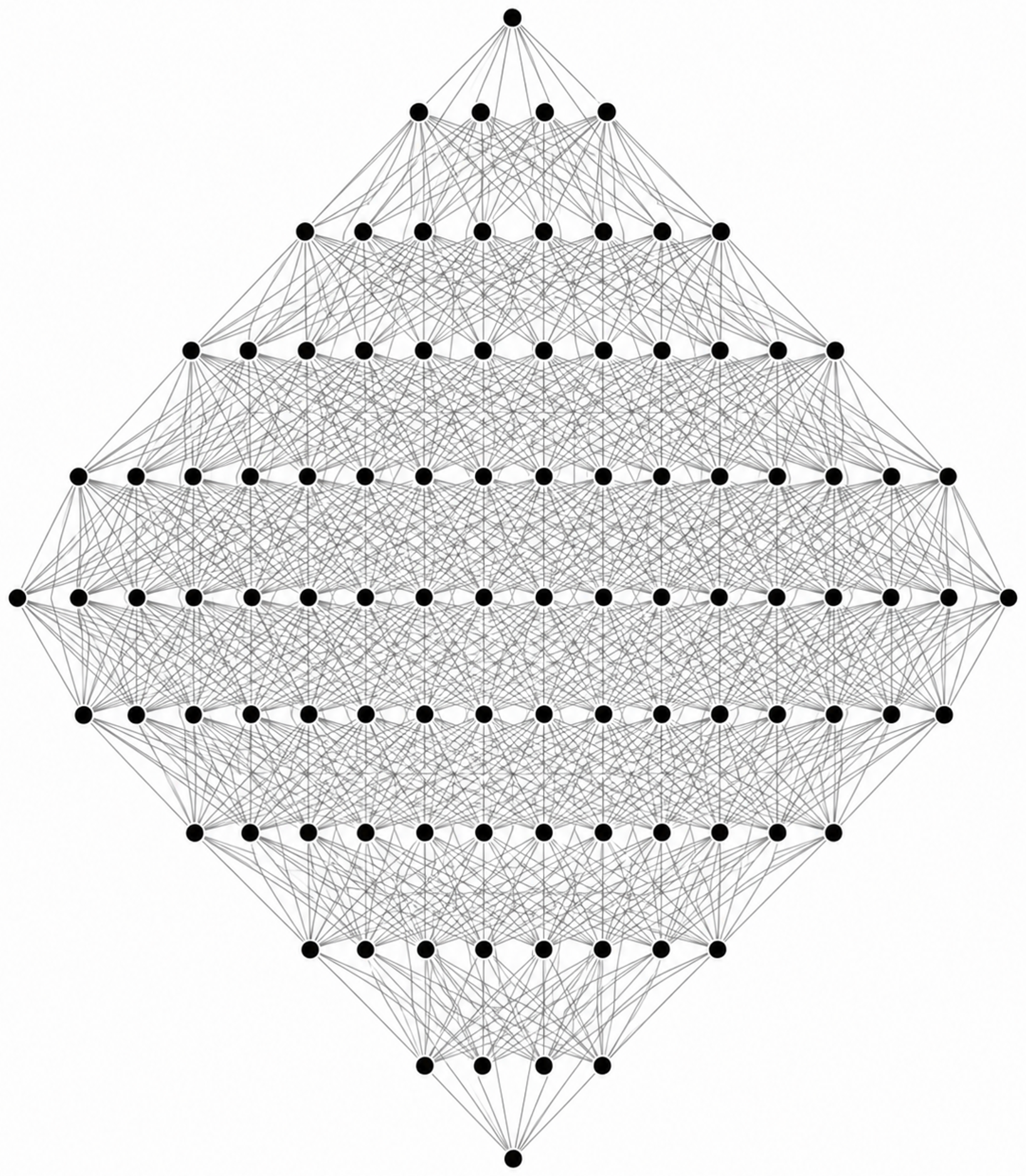}
  \includegraphics[width=0.3\columnwidth]{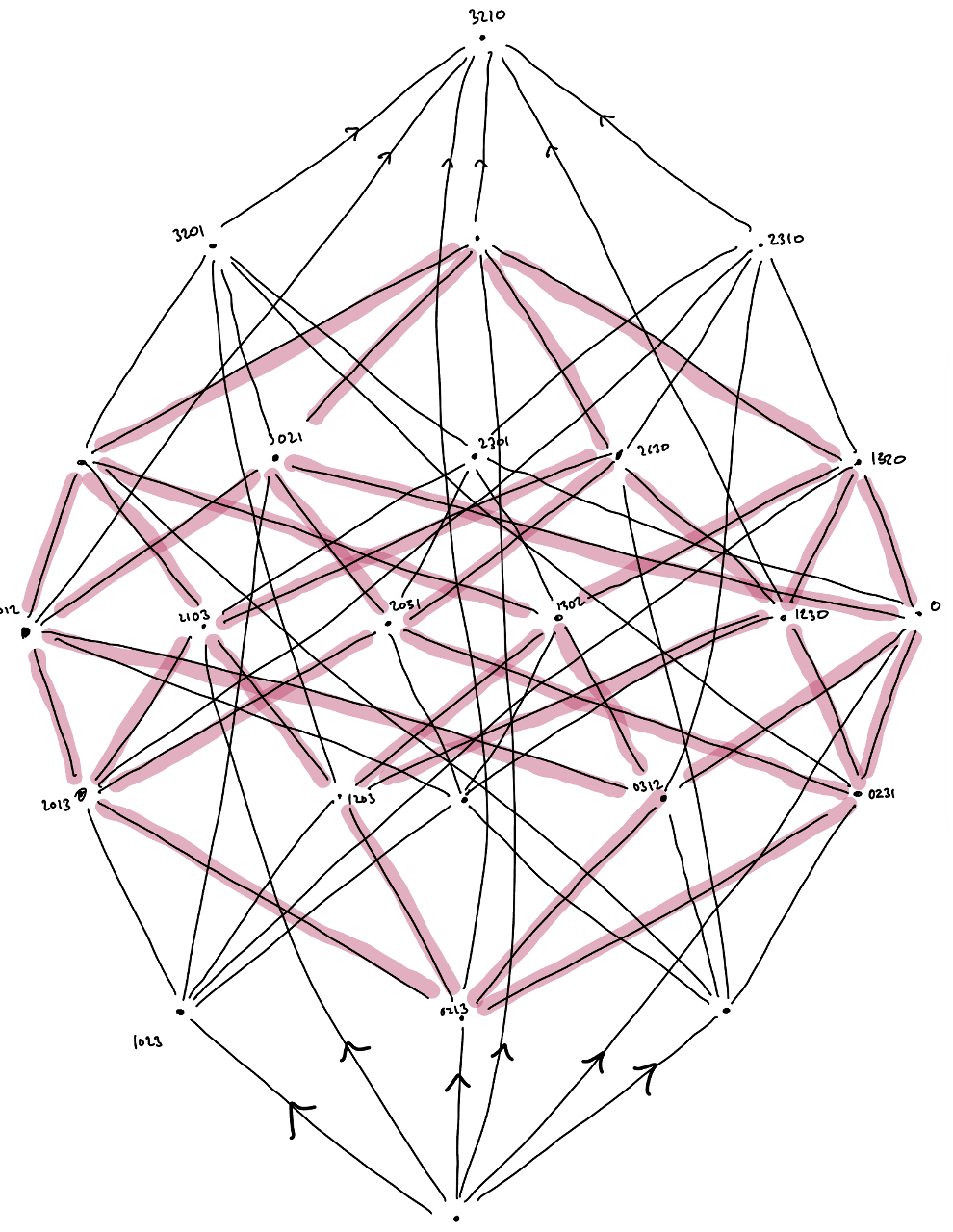}
  \caption{The Cayley graph of $S_n$ for $n=4$ (left) and $n=5$
    (middle), and a large hypercube (right).}
  \label{fig:cayley}
\end{figure}

Consider the symmetric group $S_n$ consisting of all permutations of
$1, 2, \dots, n$. This group is generated by transpositions $t_{ij} = (i,j)$,
and we consider the Cayley graph $\Gamma_n$ of $S_n$ with respect to
these generators: its vertices are all
permutations, and two vertices $x$ and $y$ are connected by an
edge if $x = t_{ij}y$. The graph $\Gamma_n$ is a highly structured graph
with $n!$ vertices (see Figure \ref{fig:cayley}).

\begin{question}
  What is the largest full subgraph of $\Gamma_n$ that is isomorphic
  to a hypercube? What is the largest dimension of such a hypercube
  that can occur? 
\end{question}

We won't motivate this question, but ask for the reader's trust that
this question is interesting (see the introduction to \cite{ELPSW} for 
motivation). Some evidence that this is an interesting question is
that small examples (where one can search exhaustively) reveal
interesting and ``unexpected'' hypercubes. (See Figure \ref{fig:cayley}
for a 4-hypercube inside the Cayley graph for $S_4$---this is an example of the kind of
subgraph we are looking for.)

On the face of it this is a very challenging search problem. A
fundamental idea espoused in \cite{romera2024mathematical} is to search for objects which
admit concise descriptions via Python programs.\footnote{In more
  mathematical terms, we are
searching for subgraphs of low Kolmogorov complexity.} This is
typically a difficult space to search over. However, LLMs give us
a sense in which two Python programs are ``similar''. For example, it
is possible to give an LLM several Python programs and ask it to
produce some other programs which are similar to the programs it has
already seen. This allows one to set up a search
amongst Python programs along the lines of genetic algorithms, which
is the basis of the \emph{FunSearch} and \emph{AlphaEvolve}
algorithms \cite{romera2024mathematical,novikov2025alphaevolve}.

One sets up an experiment by specifying what type of programs one is
searching for, and a score function that should be maximised. In our
setting, our evaluation computed the size of the hypercube for $n$
ranging from 10 to 50, and then took the average.\footnote{More precisely,
  we were optimising a related quantity called the ``$d$-invariant'',
  which in our setting turns out to be equal to the dimension of the
  hypercube, see \cite{ELPSW}.} Several promising experiments with
FunSearch (which is open source) found hypercubes of much larger
dimension than we expected. We then switched to a version of
AlphaEvolve that was internal to DeepMind. Overnight, AlphaEvolve 
found a program that was performing much better than anything we had
found earlier. Moreover, it appeared to be producing large hypercubes for
all $n$.

We then found ourselves in an interesting position: we had a program that
performs well but we didn't understand why. It was then our job to
``decompile'' the Python program into something with human
meaning. This took several days, and the result is the construction of
a natural hypercube inside $\Gamma_n$ of dimension growing like $n
\log n$, which is within a constant factor of best possible \cite{ELPSW}.

\begin{remark}
  This is one example where AlphaEvolve was successful, and I think
  illustrates an interesting application of LLMs to mathematics. In
  many other examples it does not perform as well. The
  reader is referred to \cite{AE} for a much more comprehensive study
  of the effectiveness of the AlphaEvolve algorithm on problems in
  mathematics.
\end{remark}


\section{The concern.} In July, when this lecture was delivered, we
were digesting the first solid evidence that AI systems can make
progress on major open problems (unit distance conjecture \cite{AlonEtAl2026UnitDistance}, cycle double cover
conjecture \cite{OpenAI2026CycleDoubleCover,Oum2026CycleDoubleCover}
and Jacobian conjecture counter-example \cite{Gao2026Jacobian}). Since July, we have been
flooded with major results: 
existence of non sofic groups \cite{OpenAI2026TenAdvances}, possible
complex structures on $S^6$ \cite{Alpoge2026S6,Engel2026S6} and of course, a
blowing up solution to the forced Euler equation
\cite{AlpogeBuckmaster2026Euler} and to the forced Navier-Stokes problem in fluid dynamics
\cite{OpenAI2026NavierStokes} (a Clay Millennium problem). These
breakthroughs have been accompanied by wide-spread changes in
mathematical practice: a large subset of working mathematicians now routinely use
AI systems for everything ranging from grammar checking to proving
theorems. As I wrote in the introduction, these developments have led to an
explosion of commentary from mathematicians and the broader
public. Here I would like to expand slightly on two points which were
made in July.


\subsection{Understanding and the difficulties ahead.} A simple model
for the difficulties we are currently facing is provided by Goodhart's
law: ``when a measure becomes a target, it ceases to be a good
measure''.  It is difficult to precisely pin down what our goal as mathematicians
is (and this varies from mathematician to mathematician), but an excellent approximation is given by
\emph{understanding}. 

The fundamental difficulty is that it is hard to assess
understanding. For several centuries, mathematicians have settled upon
proof as a surrogate measurement of understanding. Here is Hilbert, from his article where 
he put forward his famous list of problems \cite{Hilbert1900ICM}:
\begin{quote}
  “...rigor, which has become proverbial in mathematics, corresponds
  to a universal philosophical necessity of our understanding...”
\end{quote}
It is remarkable how closely coupled proof, rigor and understanding have
remained over the centuries. It is rare to have proof without
understanding. As a community we are generally suspicious of
purported understanding that cannot be made rigorous. Often new
understanding is accompanied by new proofs of theorems.

Over the last 50 years we have been witnessing a gradual
decoupling of proof and understanding. The canonical example of this
decoupling is the proof of the four colour theorem
\cite{AppelHaken1977a,AppelHakenKoch1977}. Here a central open problem
in topology was reduced (via considerable cleverness) to 
an exhaustive check of thousands of cases, which was completed on
computer. Here is Thurston, commenting on the controversy that
ensued (his emphasis) \cite{Thurston1994}: 
\begin{quote}
  ``I interpret the controversy as having little to do with doubt
  people had as to the veracity of the theorem or the correctness of
  the proof. Rather, it reflected a continuing desire for \emph{human
  understanding} of a proof, in addition to knowledge that the theorem
  is true.''
\end{quote}
Over the last fifty years we have been confronted with many more
such examples. A striking example (now already a decade old) is the proof of the
Pythagorean triple conjecture, whose proof is computer generated and
amounts to 200TB (trillions of lines!) \cite{HeuleKullmannMarek2016}.

As a community, we have grudgingly accepted these examples,
but not fully internalised their lesson. Namely, that increasingly
powerful systems will provide proofs without understanding. As a
community we must restructure our incentives to better reflect what we
actually value. We also have to communicate what we value to the
public and to policy makers. This will be a painful process!

\subsection{Expect more powerful systems.} Several mathematicians have
told me that they expect big AI labs to ``give up on mathematics
soon, and turn their attention elsewhere''. The prevailing sentiment
seems to be that the seismic developments of the last months will
settle down, and mathematics will return to being somewhat similar to
what it has been for the last fifty years. I strongly believe this is
wishful thinking. Rather, I think that the most rational hypothesis at the
moment is that AI systems will continue to improve.

As I tried to stress earlier, our current understanding of why models
sometimes perform so well is limited. Thus, we have very little
understanding of the mechanism behind modern AI systems, and hence poor
mental models of their future capabilities. In the absence of
theoretical understanding, the best indicators appear to be the
performance on benchmarks correlated with performance in mathematics.

\begin{figure}[tbp]
  \centering
  \includegraphics[width=0.9\columnwidth]{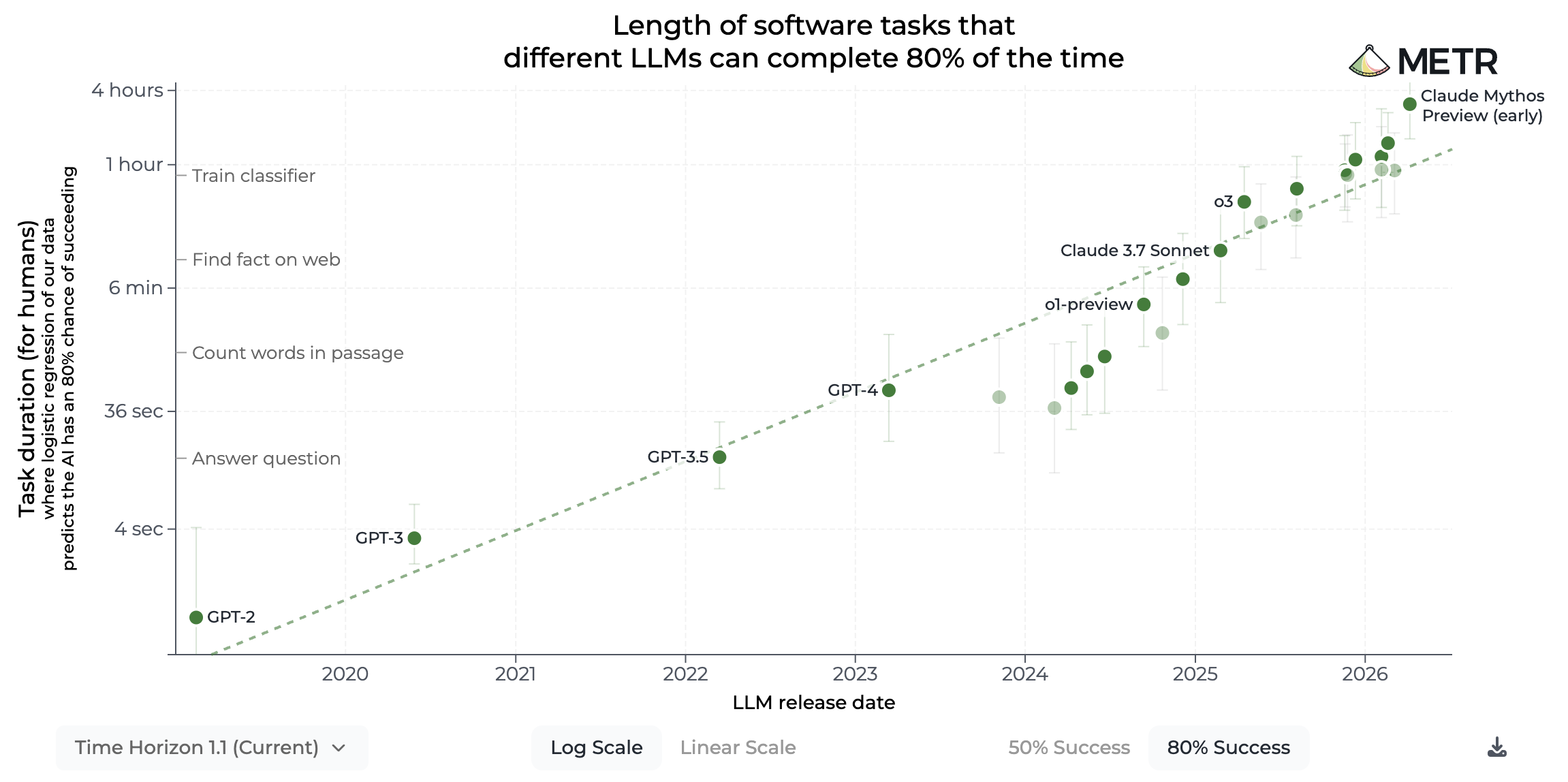}
  \caption{Time horizon of different LLMs \cite{METRTimeHorizons}.}
  \label{fig:time_horizon}
\end{figure}

Research mathematics is one of the slowest of all human processes. It
is not unusual for a researcher to spend decades working on and off on a
problem, sometimes without clear feedback as to whether they are working in the
right direction. For this reason, a critical measurement of task
difficulty appears to be the amount of time it takes a \emph{human} to
perform the task. In a famous study, METR examined the
capabilities of models over time, in terms of this measurement of task
difficulty. The conclusion seems to be that the (human) task length
that models can complete is growing exponentially (see Figure
\ref{fig:time_horizon}, noting the log scale on the $y$-axis).

\begin{figure}[tbp]
  \centering
  \includegraphics[width=0.5\columnwidth]{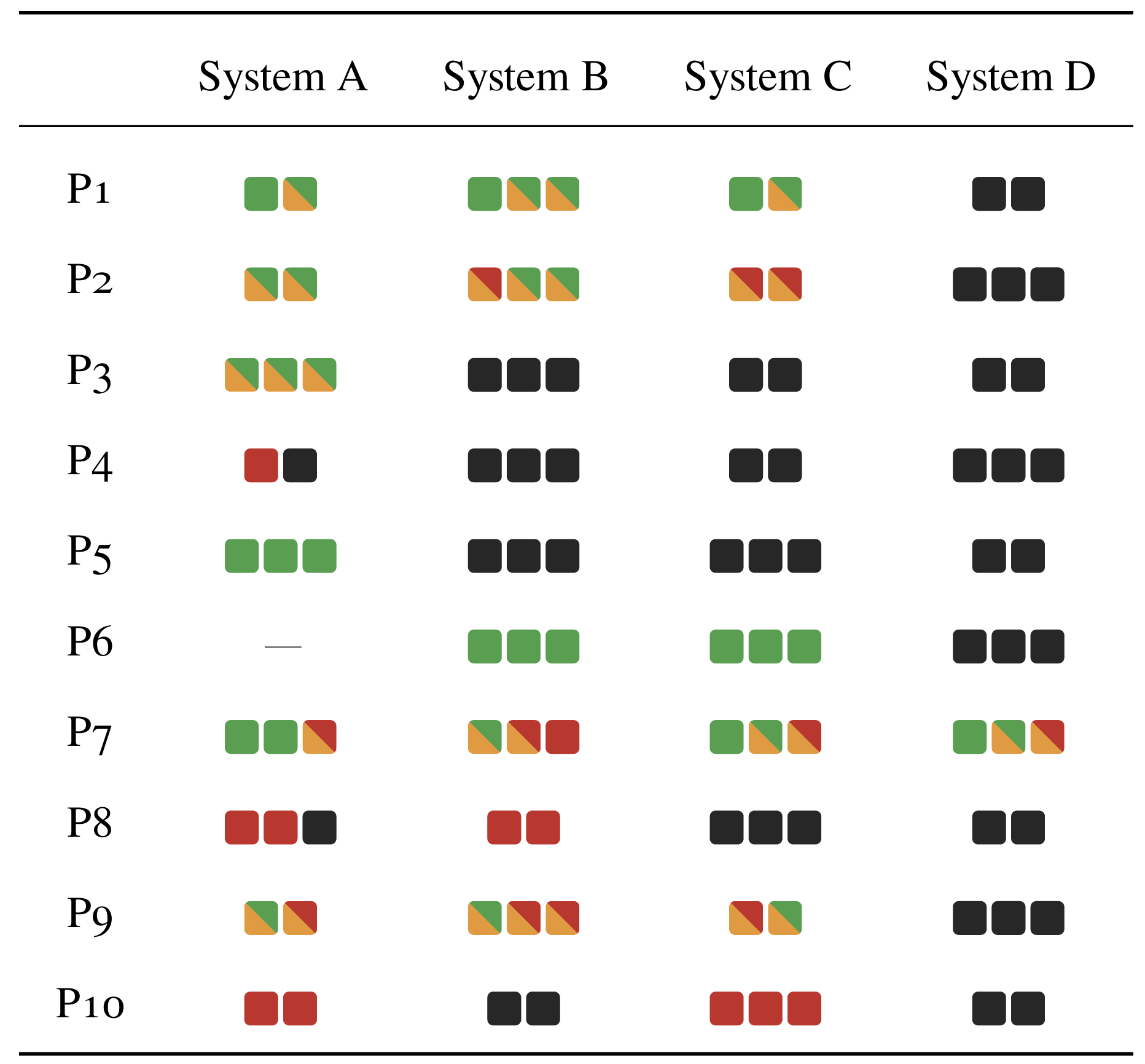}
  \caption{The results of First Proof, Second Batch \cite{AbouzaidEtAl2026FirstProofSecondBatch}.}
  \label{fig:second_batch}
\end{figure}

For research mathematicians, the most comprehensive picture is
probably provided by the First Proof
benchmark \cite{AbouzaidEtAl2026FirstProof,
  AbouzaidEtAl2026FirstProofSecondBatch}---the only problem is that it
has not been running long
enough to enable meaningful extrapolation. In this benchmark, research
mathematicians are asked to submit challenging statements from work in
progress that has not yet been published, alongside proofs and a
commentary on the perceived difficulty. Problems are then selected,
and submitted harnesses are given one shot at solving them. Solutions are then
refereed by human referees, and judged as flawless, needing minor or
major revisions, or rejected. The results of the second batch (see
Figure \ref{fig:second_batch}) give an excellent picture of LLM
capability in June 2026. The third batch is currently underway.

\section{Conclusion.} This is a time of great change in our
discipline, and in the world. I have tried to discuss these changes in
the context of the long conversation of mathematics, emphasising that
some of the fundamental ideas underlying modern AI systems arose within
this long conversation. I believe that a greater awareness of the many
mysteries raised by AI systems will benefit our discipline. Major changes lie
ahead, but I strongly believe that human mathematics will
survive and flourish. Let me finish with another quote from Thurston \cite{Thurston2010Mathematician}:
\begin{quote}
  ``\dots mathematics only exists in a living community of mathematicians that spreads understanding and breathes life into ideas both old and new.''
\end{quote}

\section*{AI use statement.} No AI was used for writing. ChatGPT
was used for references and for proof-reading.

\section*{Acknowledgments.}
I would like to thank
Pramod Achar, Alex Davies, Keeley Hoek, Tom Goertzen, Nalini Joshi, Finn Klein, Daniel Murfet,
Simon Riche, Connor Simpson, John Voight, Oded Yacobi and Tamar
Ziegler for very useful feedback on the lecture and these notes.

\bibliographystyle{siamplain}
\bibliography{references}
\end{document}